\documentclass[9pt]{article}
\usepackage{color}
\usepackage{amsmath}
\usepackage{amsthm}
\usepackage{wrapfig}
\usepackage{verbatim}
\usepackage{amsfonts}
\usepackage{dsfont}

\usepackage{tikz}
\usepackage{pgf}
\usepackage{bbm}
\usepackage{centernot}
\usetikzlibrary{arrows}
\usetikzlibrary{decorations.pathmorphing}

\begin{document}
{\sc Arman Margaryan\par
 Ners\'es Aramian}\vspace{0.5cm}
\begin{center}
{\large Compass Constructions}
\end{center}\vspace{0.5cm}\par
 The purpose of this short manuscript is to show that all point constructions that can be done via ruler and compass, can also be done with compass exclusively. If we are using compass and ruler the way we construct new points is by first constructing circles or lines and then considering their intersection. This hints us a way to approach the problem, which is possibly the most straightforward approach. However, the execution of it is rather interesting and worth the effort. The approach that we are talking about is to show that the intersection of circles and lines, lines and lines can be constructed using only compass. The intersection of two circle comes in as trivial given, since all we need to do is construct circles in this situation. Our task can be formulated even more precisely. Suppose we are given four points $A$, $B$, $C$ and $D$. Our task is to show that the intersection of lines $AB$ and $CD$ can be constructed using only compass, and that the intersection of $AB$ and the circle centered $C$ and passing through $D$ can be also constructed using compass.\par
Rather surprisingly the first hard task is to construct the midpoint of a segment. As a matter of fact, the construction itself is not difficult at all and it is shown on the next page. The difficult part is coming with the construction. Here we will take a step by step approach, which will enhance our understanding of the thinking process behind the construction. We think it will beneficial to the reader to go through these ideas. However, the reader may skip the discussion, and go to the sketch of the construction on the next page. It will not disrupt the continuity of the flow.\par
Suppose we are given an interval on the plane. As it is common in mathematics we will identify the plane with the set of complex numbers, $\mathbb{C}$, where one of the endpoints of the interval is identified with $0$ and the other with $1$. This is just temporary and we are going to resume our regular planar geometry language shortly after. In this formulation, our task is simply to construct $\frac{1}{2}$ using only compass. To make the matters even more pedantic, let us denote the set of compass-constructable points with $\mathfrak{C}\subset \mathbb{C}$. Thus, our task is to show that $\frac{1}{2}\in\mathfrak{C}$.\par

Our first claim does not exploit the fact that we are working with the complex numbers; it is rather a geometric construct. Suppose that $a,b\in \mathfrak{C}$, and $c$ is a point in plane, such that the resulting triangle is equilateral, then $c\in \mathfrak{C}$. This claim is rather obvious and should be easy to check. The crucial use of this claim comes in
\begin{wrapfigure}{r}{0.3\textwidth}
\begin{center}
\begin{tikzpicture}
\draw[red] (0,0) -- (0.5,0.866) -- (1,0) -- (2,0) -- (1.5,0.866) -- (1,0)
(1.5,0.866) -- (0.5,0.866) -- (-0.5,0.866) -- (0,0) -- (-1,0) -- (-0.5,0.866);
\draw[fill] (0,0) circle (0.025) node[below] {\footnotesize $0$};
\draw[fill] (1,0) circle (0.025) node[below] {\footnotesize $1$};
\draw[fill,red] (-1,0) circle (0.025) node[below left] {\footnotesize $-1$};
\draw[fill,red] (0.5,0.866) circle (0.025) node [above right] {\footnotesize $\omega$};
\draw[fill,red] (1.5,0.866) circle (0.025) node [above right] {\footnotesize $\omega+1$};
\draw[fill,red] (-0.5,0.866) circle (0.025) node [above left] {\footnotesize $\omega-1$};
\draw[fill,red] (2,0) circle (0.025) node [below right] {\footnotesize $2$};
\draw (0,0) -- (1,0);
\end{tikzpicture}
\end{center}
\end{wrapfigure} 
proving that $-1,2\in \mathfrak{C}$. The figure on the right demonstrates the construction using the equilateral triangles. This result will enable us to show the following key lemma.\par
\underline{\sc Lemma.} {$\mathfrak{C}$ is a subring of $\mathbb{C}$, which is closed under conjugation.}\par
\underline{\sc Proof.} The clear starting point is to show that $0,1\in \mathfrak{C}$; however, that is the case, since those are our starting points. If we show now that $\mathfrak{C}$ is closed under addition, multiplication and conjugation then we will be done. Indeed, the only missing part will be the existence of additive inverses; however, that is fixed by the fact that $-1\in \mathfrak{C}$ and multiplication by it, additively inverts the number.\par
Conjugation is the easiest of all. The picture on the right shows quite clearly how to conjugate a point in 
\begin{wrapfigure}{r}{0.3\textwidth}
\begin{center}
\begin{tikzpicture}
\draw[red] (0,0) ++ (60:1.732) arc (60:-60:1.732);
\draw[red] (1,0) ++ (80:1) arc (80:-80:1);
\draw[fill] (0,0) circle (0.025) node[below] {\footnotesize $0$};
\draw[fill] (1,0) circle (0.025) node[below] {\footnotesize $1$};
\draw[fill] (1.5,0.866) circle (0.025) node[above right] {\footnotesize $a$};
\draw[fill,red] (1.5,-0.866) circle (0.025) node[below right] {\footnotesize $a^*$};
\end{tikzpicture}
\end{center}
\end{wrapfigure} 
$\mathfrak{C}$. For a given $a\in \mathfrak{C}$, simply intersect the circles centered at $0$ and $1$, and passing through $a$. The set of intersection points is $\{a,a^*\}$ -- quite easy to check.\par
Multiplication is also relatively easy to construct.
Suppose that in addition to $a$, $b$ is also compass-constructable. We can exactly repeat the construction of $b$, but this time around treating $0$ and $a$ as our starting points (as opposed to $0$ and $1$). The resulting point, $c$, will have the following two obvious properties. The change of initial points scales the lengths by a factor of $|a|$. Therefore, we conclude $|c|=|b||a|$. The second property is that the phase angle difference between $c$ and $a$ is the same as the phase angle of $b$. The conclusion of these two statements is that $c=ab$.\par
To show that $\mathfrak{C}$  is closed under addition we use the trick of changing the initial points and mimicking the previous construction. However, this time we do it twice. Begin by mimicking the construction of $a$ with respect to starting points $1$ and $2$. This clearly results in construction of point $a+1$. Now mimic the construction of $b$ with respect to initial points $a$ and $a+1$. This does the proof, since the construction creates the desired point $a+b$. $\qed$
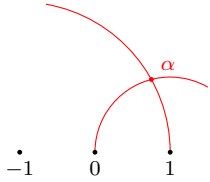
\par

\begin{wrapfigure}{l}{0.3\textwidth}
\begin{center}
\begin{tikzpicture}
\draw[red] (1,0) arc (0:80:2);
\draw[red] (0,0) arc (-180:-300:1);
\draw[fill,red] (0.75,0.968) circle (0.025) node[above right] {\footnotesize $\alpha$};
\draw[fill] (0,0) circle (0.025) node[below] {\footnotesize $0$};
\draw[fill] (1,0) circle (0.025) node[below] {\footnotesize $1$};
\draw[fill] (-1,0) circle (0.025) node[below] {\footnotesize $-1$};
\end{tikzpicture}
\end{center}
\end{wrapfigure} 
We are almost there, but only using the lemma above we cannot show the existence of construction. We need an extra step, an extra point constructed. Consider the figure on the left. The point $\alpha$ is  an intersection of two circles centered at $-1$ and $1$ with radii $2$ and $1$ respectively. We let $\alpha$ be $\frac{3+i\sqrt{15}}{4}$, which is one of the points of intersection. According to the lemma $|\alpha|^2=\frac{3}{2}\in \mathfrak{C}\Longrightarrow \frac{1}{2}\in\mathfrak{C}$. If one pays attention, we did not need the lemma to its full power to prove our original claim; however, it is an interesting fact on its own. We thought it may be interesting.\par
\begin{wrapfigure}{r}{0.3\textwidth}
\begin{center}
\begin{tikzpicture}[scale=0.8]
\draw[red] (2,0) circle (1);
\draw[red] (0,0) circle (2);
\draw[red] (1.75,0.968) circle (1);
\draw[red] (1.75,-0.968) circle (1);
\draw[red] (1,0) circle (1);
\draw[red] (1.5,0.866) circle (1);
\draw[red] (0.5,0.866) circle (1);
\draw[fill] (1,0) circle (0.025) node[below] {\footnotesize};
\draw[fill] (2,0) circle (0.025) node[below left] {\footnotesize};
\draw[fill,red] (1.5,0) circle (0.025) node[below left] {\footnotesize};
\end{tikzpicture}
\end{center}
\end{wrapfigure} 
The construction follows straight from our discussion and is depicted below. We will recommend the reader to examine as an exercise how each circle is constructed and convince himself that the construction is valid.\par
Now we can abandon complex numbers and go to regular planar geometry. It is rather surprising how useful the fact that we have just proved is. One trivial consequence is that for any two distinct points we can construct the circle that has the segment between those two points as its diameter. This enables us to construct the bases for the perpendiculars. Specifically, if we are given point three points, $A$, $B$, $C$, such that $A\neq B$, then it is possible to construct the base, $H$, of the perpendicular line to $AB$ passing through $C$. The figure below demonstrates how to construct $H$: simply construct two circles one with diameter $AC$ and $BC$, and $H$ will be the intersection point on the line $AB$.\par
Now suppose we are given a circle $\Omega$ centered at point $O$, and a point $P$ outside of $\Omega$. What the figure below on the right demonstrates is the construction of inversion of $P$ with respect to $\Omega$. As the figure shows we construct the circle with diameter $OP$. We label the intersection points of the circle with $\Omega$ by $M$ and $N$. The fact that $M$, $N$ exist and distinct comes from the fact that $P$ is not in $\Omega$. Then the base of the perpendicular from $O$ to $MN$, which is denoted by $I$, will be our desired inversion point.\par
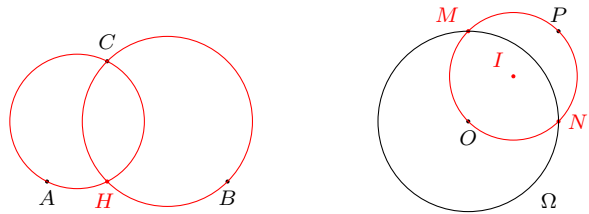
\begin{center}
\begin{tikzpicture}[scale=0.8]
\draw[fill] (0,0) circle (0.025) node[below] {\footnotesize $A$};
\draw[fill] (1,2) circle (0.025) ++ (0,0.05) node[above] {\footnotesize $C$};
\draw[fill] (3,0) circle (0.025) node[below] {\footnotesize $B$};
\draw[fill,red] (1,0) circle (0.025) ++ (-0.05,-0.05) node[below] {\footnotesize $H$};
\draw[red] (0.5,1) circle (1.118);
\draw[red] (2,1) circle (1.414);
\begin{scope}[xshift=7cm, yshift=1cm]
\draw[fill] (0,0) circle (0.025) node[below] {\footnotesize $O$};
\draw[fill] (1.5,1.5) circle (0.025) node[above] {\footnotesize $P$};
\draw[fill,red] (1.5,0) circle (0.025) node[right] {\footnotesize $N$};
\draw[fill,red] (0,1.5) circle (0.025) node[above left] {\footnotesize $M$};
\draw[fill,red] (0.75,0.75) circle (0.025) node[above left] {\footnotesize $I$};
\draw (1.05,-1.05) node[below right] {\footnotesize $\Omega$};
\draw[red] (0.75,0.75) circle (1.061);
\draw (0,0) circle (1.5);
\end{scope}
\end{tikzpicture}
\end{center}\par
What about the rest of the points? Is it possible to invert them? The answer is obviously positive for the points on $\Omega$. For the points inside $\Omega$, the answer is still positive; however, the construction is not as ``universal'' as it was the case with the exterior point. More specifically, the construction is actually point dependent. Here is how one must proceed. $P$ must be different from $O$ otherwise the inversion is not defined. From the lemma, we can see that we can construct point $Q$, such that $P$ is on the segment $OQ$ and that the ratio of length of segment $QO$ to the length of segment $OP$ is a positive integer. Let the integral ratio be so large that $Q$ is outside of $\Omega$. Construct the inversion of $Q$ and call it $J$. The inversion of $P$, which again will be called $I$, will be the point that satisfies the properties: $J$ belongs to the segment $OI$ and the ratio of the length of the segment $OI$ to the length of $OJ$ is the same positive integer as before. Overall, we conclude that the inversions can be constructed.\par
Now we can tackle the original questions proposed in the beginning of the manuscript dealing with the intersections of lines and intersections of lines and circles. We start off with the intersections of lines. Suppose we are given the points $A$, $B$, $C$ and $D$, such that they are not colinear and $AB \centernot\parallel CD$. We want to find the intersection of $AB$ and $CD$, which we will call $S$. We hope that no one will have an objection against us stating that it is possible to construct a point $P$, such that $P\not\in AB\bigcup CD$, and a circle $\Gamma$ centered at $P$. 
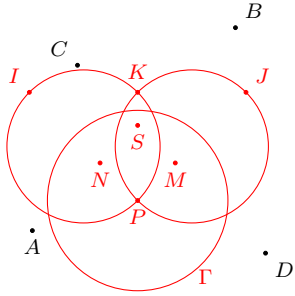
\begin{wrapfigure}{l}{0.3\textwidth}
\begin{center}
\begin{tikzpicture}
\draw[fill] (-0.4,-0.4) circle (0.025) node[below] {\footnotesize $A$};
\draw[fill] (2.3,2.3) circle (0.025) node[above right] {\footnotesize $B$};
\draw[fill] (0.2,1.8) circle (0.025) node[above left] {\footnotesize $C$};
\draw[fill] (2.7,-0.7) circle (0.025) node[below right] {\footnotesize $D$};
\draw[fill,red] (1,0) circle (0.025) node[below] {\footnotesize $P$};
\draw[fill,red] (0.5,0.5) circle (0.025) node[below] {\footnotesize $N$};
\draw[fill,red] (-0.44,1.44) circle (0.025) node[above left] {\footnotesize $I$};
\draw[fill,red] (1.5,0.5) circle (0.025) node[below] {\footnotesize $M$};
\draw[fill,red] (2.44,1.44) circle (0.025) node[above right] {\footnotesize $J$};
\draw[fill,red] (1,1.44) circle (0.025) ++ (0,0.05)node[above] {\footnotesize $K$};
\draw[fill,red] (1,1) circle (0.025) node[below] {\footnotesize $S$};
\draw[red] (1.9,-0.8) node[below] {\footnotesize $\Gamma$};
\draw[red] (1,0) circle (1.2);
\draw[red] (0.28,0.72) circle (0.707*1.44);
\draw[red] (1.72,0.72) circle (0.707*1.44);

\end{tikzpicture}
\end{center}
\end{wrapfigure}
First begin by constructing the bases of perpendiculars from $P$ to $AB$ and $CD$. Call them $N$ and $M$ respectively. Invert $N$ to $I$ relative to $\Gamma$. Do the same for $M$ to obtain point $J$. Now draw two new circles that have $IP$ and $JP$ as diameters. Let $K$ denote the intersection point of the circles that is distinct from $P$. Now invert $K$ with respect to $\Gamma$ and what we obtain is, in fact, $S$. The reader should verify that all the steps that we have done were correct and that the endresult is what we were looking for.\par
The intersection of a circle and a line actually takes more effort. This is mainly due to the fact that there are two cases to consider. Consider the points $A$, $B$ and a circle $\Omega$ with center $O$, such that $AB$ intersects $\Omega$. For the first case assume that $O\not\in AB$. 
We apply again the inversion trick and obtain the intersection points of $AB$ with $\Omega$. This one is actually  much shorter. Construct the base of the perpendicular from $O$ to $AB$. Call it $H$. From our assumption it follows that $H$ is distinct from $O$; therefore, we can invert it with respect to $\Omega$. Call this inverted point $K$. Draw another circle with diameter $OK$. The intersection point(s) of this circle with $\Omega$ are actually the points of intersection of $AB$ with $\Omega$. Again it is advisible to the reader to validity of the claim on his own.
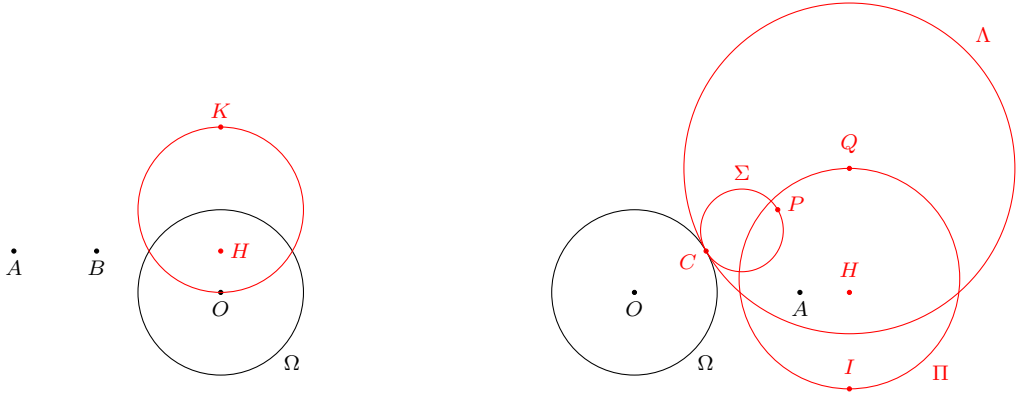
\begin{center}
\begin{tikzpicture}[scale=1.1]
\draw[fill] (-2.5,0.5) circle (0.025) node[below] {\footnotesize $A$};
\draw[fill] (-1.5,0.5) circle (0.025) node[below] {\footnotesize $B$};
\draw[fill] (0,0) circle (0.025) node[below] {\footnotesize $O$};
\draw[fill,red] (0,0.5) circle (0.025) node[right] {\footnotesize $H$};
\draw[fill,red] (0,2) circle (0.025) node[above] {\footnotesize $K$};
\draw (0.65, -0.65) node[below right] {\footnotesize $\Omega$};
\draw (0,0) circle (1);
\draw[red] (0,1) circle (1);
\begin{scope}[xshift=5cm]
\draw[fill] (2,0) circle (0.025) node[below] {\footnotesize $A$};
\draw (0.65, -0.65) node[below right] {\footnotesize $\Omega$};
\draw[fill] (0,0) circle (0.025) node[below] {\footnotesize $O$};
\draw[fill,red] (0.866,0.5) circle (0.025) ++ (0,0.07) node[below left] {\footnotesize $C$};
\draw[fill,red] (0.866*2,1) circle (0.025) ++ (0,0.07) node[right] {\footnotesize $P$};
\draw[fill,red] (0.866*3,1.5) circle (0.025) ++ (0,0.07) node[above] {\footnotesize $Q$};
\draw[fill,red] (0.866*3,1.5-4/1.5) circle (0.025) ++ (0,0.07) node[above] {\footnotesize $I$};
\draw[fill,red] (0.866*3,0) circle (0.025) ++ (0,0.07) node[above] {\footnotesize $H$};
\draw[fill,red] (0.866*3,1.5) ++ (1.41,1.41) node[above right] {\footnotesize $\Lambda$};
\draw[fill,red] (0.866*3/2,1.5/2) ++ (0,0.5) node[above] {\footnotesize $\Sigma$};
\draw[fill,red] (0.866*3,1.5-4/1.5) ++ (1.1,0)  node[above] {\footnotesize $\Pi$};
\draw (0,0) circle (1);
\draw[red] (0.866*3, 1.5) circle (2);
\draw[red] (0.866*3/2, 1.5/2) circle (0.5);
\draw[red] (0.866*3, 1.5-2/1.5) circle (2/1.5);
\end{scope}
\end{tikzpicture}
\end{center}\par
Now we are going to tackle the case where $O\in AB$, which is surprisingly quite tricky, at least as far as we know. As a matter of fact, it can be used to solve the previous problem as well. We would like to thank Daniel Briggs for help on this last part. For this case $B$ is redundant; therefore, we are going to ignore it. We want to find the intersection points of $OA$ and $\Omega$. We think that nobody will complain if we state that starting with $O$ and $A$ we can construct a point $C\in \Omega$, such that $C\not\in OA$. Then construct $P$ and $Q$, so that $OC=CP=PQ$. Draw a circle centered at $Q$ and passing through $C$. Call that circle $\Lambda$. Then consider $H$, the base of the perpendicular from $Q$ to $OA$. Invert it with respect to $\Lambda$, and call the inversion $I$. Let $\Pi$ be the circle with diameter $QI$, and $\Sigma$ be the circle with diameter $CP$. One then can easily convince himself that if we invert the two points of $\Sigma\cap \Pi$ with respect to $\Lambda$, then we will obtain the intersection points of $OA$ and $\Omega$.\par
This accomplishes the goal of the manuscript. As you may have noticed the manuscript does not contain many proofs, and to be frank, they are mostly fairly straightforward. However, coming up with these constructions was not an easy task and required some significant effort. Hopefully, you enjoyed reading this manuscript through.
\end{document}